\documentclass[12pt]{article}
\usepackage{latexsym, amsmath, amssymb, amsfonts}
\title{A Note on a Modified Catuskoti}
\author{Tetu Makino\footnote{Department of Applied Mathematics, Yamaguchi University, Ube 755-8611, Japan. 
E-mail makino@yamaguchi-u.ac.jp}}
\date{\today}
\begin{document}
\maketitle

\begin{abstract}
The `catuskoti' or tetralemma in Buddhist logic is a problematic subject from the modern logical point of view. Recently a many-valued paraconsistent logic was proposed in order to formalize catuskoti adequately by G. Priest. On the other hand a slight modification of the formalization of catuskoti seems to allow an appropriate interpretation in the framework of the classical propositional calculus in the mathematical logic developed by Russell-Whitehead and Hilbert-Ackermann.


\end{abstract} 

\section{Introduction}

The catu\d{s}ko\d{t}i or tetralemma in Buddhist logic is the tuple of the four alternatives

\begin{equation}
A,\quad \mbox{not}\  A, \quad \mbox{both}\  A\  \mbox{and not}\  A,
\quad \mbox{neither}\  A\  \mbox{nor not}\  A,
\end{equation}
where $A$ is a proposition or a predicate. From the modern logical point of view it was investigated by K. N. Jayatilleke \cite{Jayatilleke} (1967), D. S. Ruegg \cite{Ruegg} (1977), R. D. Gunaratne \cite{Gunaratne} (1986), J. Westerhoff \cite{Westerhoff} (2006) and others. If we insist the classical two-valued logic, this tuple may be formulated as
\begin{equation}
A, \quad \neg A, \quad A\wedge (\neg A), \quad \neg(A\vee(\neg A)).
\end{equation}
Here and hereafter we adopt the propositional calculus developed in 
\cite{Sider}. The definitions of symbols and terminology and axioms follow it
except for using the symbol $\neg$ instead of $\sim$. However the fourth alternative $\neg(A\vee(\neg A))$ of (2) is equivalent to $(\neg A)\wedge A$ and to the third alternative if we are confined to the classical logic. Therefore G. Priest \cite{Priest} (2010) introduced a formal logical machinery which may be more appropriate than the classical one. An adequate formalization of the catu\d{s}ko\d{t}i requires a four-valued logic. See \cite{Priest}. 

In this note, we introduce a modification of the catu\d{s}ko\d{t}i, and give an its interpretation in the framework of the classical two-valued logic in the form of the propositional calculus developed in \cite[Chapter 2]{Sider}. 

\section{Modification and interpretation}

Let us modify the tuple (2) replacing the third alternative $A\wedge(\neg A)$ by $A\vee(\neg A)$. That is, we consider the tuple of formulas
\begin{equation}
A,\quad \neg A, \quad A\vee(\neg A),
\quad \neg(A\vee(\neg A)).
\end{equation}
Here $A$ is a {\bf formula} of a propositional calculus $\Sigma$ we consider.
Of course the third formula is identically true (a tautology, $\curlyvee$), and the fourth formula is identically false (a contradiction,
$\curlywedge$) in the usual semantics of the classical propositional calculus.

The characterization of this tuple is as follows. Let us denote by $L_0$ the set of all {\bf sentence letters} of $\Sigma$ and by $L$ the set of all formulas of $\Sigma$. A mapping $V$ from $L_0$ into $\{0,1\}$ is called a {\bf valuation}. If a valuation $V$ is given, it can be uniquely extended to a mapping from $L$ into $\{0,1\}$ by dint of
the usual truth value tables. (0 stands for `false', and 1 stands for `true'.) Let us denote this unique extension by the same letter $V$ in abbreviation.

Now we assume that there are a valuation $V_0$ such that $V_0(A)=0$ and a valuation $V_1$ such that $V_1(A)=1$. In such a case we shall say that the formula $A$ is {\bf generic} in the calculus $\Sigma$. Actually it is the case if $A$ is one of the sentence letters of $\Sigma$, that is, a member of $L_0$. It is not the case if $A$ is a tautology or a contradiction.

Then for any formula $P$ the pair of truth values $(V_0(P), V_1(P))$
should be either $(0,1), (1,0), (0,0)$ or $(1,1)$. Therefore the set $L$ of all formulas of $\Sigma$ is divided into the following four subsets:
\begin{align*}
L_1&=\{ P\in L | (V_0(P),V_1(P))=(0,1) \}, \\   
L_2&=\{ P\in L | (V_0(P),V_1(P))=(1,0) \}, \\   
L_3&=\{ P\in L | (V_0(P),V_1(P))=(1,1) \}, \\  
L_4&=\{ P\in L | (V_0(P),V_1(P))=(0,0) \}.
\end{align*}
Moreover $A, \neg A, A\vee(\neg A), \neg(A\vee(\neg A))$ are representative formulas of \\ $L_1,L_2,L_3,L_4$, respectively. This is the situation of the modified tuple (3).

Hence the interpretation of this modified catu\d{s}ko\d{t}i is: {\it If 
somebody denies all these alternatives, then he/she intends to mean by abbreviation using the representatives that he/she declares that the ultimate truth or the reality cannot be described by any formula of any propositional calculus in which $A$ is a generic formula. Particularly any propositional calculus for which $A$ is a sentence letter doesn't work to describe the reality.}\\

{\bf Note. }\    If a formula $P$ is a contradiction $\curlywedge$, like $A\wedge(\neg A)$, then\\
$(V_0(P), V_1(P))=(0,0)$. But the inverse is not true, that is,  $P$ is not necessarily a contradiction when $(V_0(P), V_1(P))=(0,0)$. If $Q$ is a tautology $\curlyvee$, then $(V_0(P), V_1(P))=(1,1)$, but the inverse is not true. In fact, as an example, let us consider the case in which $A$ is a sentence letter. Let $B$ be another sentence letter such that $A\not= B$. Then there are valuations $V_0, V_1$ such that
$V_0(A)=V_0(B)=1$ and $V_1(A)=V_1(B)=0$.
Put $P=A\wedge(\neg B)$ and $Q=A\vee(\neg B)$. Then we have $(V_0(P), V_1(P))=(0,0)$ and
$(V_0(Q), V_1(Q))=(1,1)$.
But, since there is a valuation $V_2$ such that $V_2(A)=1, V_2(B)=0$ for which $V_2(P)=1$ so that $P$ is not a contradiction. Since there is a valuation $V_3$ such that 
$V_3(A)=0, V_3(B)=1$ for which
$V_3(Q)=0$ so that $Q$ is not a tautology.

The idea to consider the pairs of the truth values is related to the philosophical point of view of G. Priest \cite{Priest}. For the details see the Appendix.     

\section{Tathagata after the death}

As the first example we take a passage from `Khema Sutta' (SN44.1) \cite[p. 1381]{SN} in the `Samyutta Nikaya':
\begin{quote}
Q1: How is it, revered lady, does the Tathagata exist after death?

A1: Great king, the Blessed One has not declared this.

Q2: Then, revered lady, does the Tathagata not exist after death?

A2: Great king, the Blessed One has not declared this either.

Q3: Then does the Tathagata both exist and not exist after death?

A3: Great king, the Blessed One has not declared this. 

Q4: Well then, does the Tathagata neither exist nor not exist after death?

A4: Great king, the Blessed One has not declared this either.

Q5: Now, what is the cause and reason, why that has not been declared by the Blessed One?

A5: That form by which one describing the Tathagata might describe him has been abandoned by the Tathagata. The Tathagata is liberated from reckoning in terms form; he is deep, immeasurable, hard to fathom like the great ocean. That feeling by which one describing....That perception by which one describing...,... 

\end{quote}
Therefore our modification is to replace Q3  by \\

Q3':  Then does the Tathagata either exist or not exist after death? Can both be allowed? \\

According our interpretation we can say that A5 explains that why the reality is beyond the set of all formulas of any calculus in which ``the Tathagata exists after death" is formalized by a generic formula.   

Now the corresponding original Pali text reads
\begin{quote}
Q1: Kinnu kho ayye hoti tath\={a}gato parammara\d{n}\={a}ti.

A1: Avy\={a}kata\d{m} eta\d{m} mah\={a}r\={a}ja bhagavat\={a}:``hoti tath\={a}gato parammara\d{n}a"ti.

Q2: Kimpanayyo na hoti tath\={a}gato parammara\d{n}\={a}ti.

A2: Etampi kho mah\={a}r\={a}ja avy\={a}kata\d{m} bhagavat\={a}: ``
na hoti tath\={a}gato parammara\d{n}\={a}"ti.

Q3: Kinnu kho ayye, hoti ca na ca hoti tath\={a}gato parammara\d{n}\={a}ti.

A3: Avy\={a}kata\d{m} kho eta\d{m} mah\={a}r\={a}ja bhagavat\={a}:`` hoti ca na ca hoti tath\={a}gato parammara\d{n}\={a}"ti.

Q4: Kimpanayye, neva hoti na na hoti tath\={a}gato parammara\d{n}\={a}ti.

A4: Etampi kho mah\={a}r\={a}ja avy\={a}kata\d{m} bhagavat\={a}:
``neva hoti na na hoti tath\={a}gato parammara\d{n}\={a}"ti.

\end{quote}

Thus the tuple of the four alternatives in the Pali is
\begin{equation}
A, \quad \mbox{na}\  A,\quad A\  \mbox{ca na ca}\  A,
\quad \mbox{neva}\  A\  \mbox{na na}\  A.
\end{equation}

According to our modification the third alternative
`$A$ ca na ca $A$' is formalized by $A\vee(\neg A)$.

The essential point lies on the difference between $A\vee(\neg A)$ and $A\wedge(\neg A)$.We wonder whether the formalization as 
$A\wedge(\neg A)$ of the translation `both $A$ and not $A$' is inevitable or not.

Let us look at Chinese translations. We could not find a Chinese translation of this `Khema Sutta', but Chinese translations
by Gunabhadra of  other many suttas in the group `Avyakata Samyutta' of `Samyutta Nikaya' can be found in the Chinese `Z\'{a} \={E}h\'{a}n-j\={i}ng '. See `The Taisho Tripitaka' 99-958, -959, -960, -962. The Chinese translation of the four alternatives are:\\

(R\'{u}l\'{a}i) y\v{o}u h\`{o}us\v{i}, w\'{u} h\`{o}us\v{i}, y\v{o}u w\'{u} h\`{o}us\v{i}, f\={e}i y\v{o}u f\={e}i w\'{u} h\`{o}us\v{i}\\
or \\

 y\v{o}u h\`{o}us\v{i}, w\'{u} h\`{o}us\v{i}, y\v{o}u w\'{u} h\`{o}us\v{i}, f\={e}i y\v{o}u h\`{o}us\v{i} f\={e}i w\'{u} h\`{o}us\v{i}.\\

Thus the tuple of the four alternatives in the Chinese is
\begin{equation}
A, \quad \neg A, \quad A\  \neg A, \quad \mbox{f\={e}i}\ A\  \mbox{f\={e}i}\ \neg A.
\end{equation} 

Since `y\v{o}u' = to have, `w\'{u}'=to lack, `f\={e}i' =not, `h\`{o}us\v{i}' =after death, the conjunctions `and' ,`or' do not appear explicitly.  In usual conversations,  ``Y\v{o}u w\'{u}" (or ``Y\v{o}u m\'{e}iyou" colloquially) does not mean ``One has {\bf and} lacks", but means ``(Do you) have {\bf or} don't have?".  In the same way ``H\v{a}o buh\v{a}o", not meaning ``It's good and bad", means ``Is it good or not?", or `` How do you think?", where `h\v{a}o'=good and `b\`{u}'=not. Therefore the putting side by side without conjunction `$A\ \neg A$' should be
interpreted as $A\vee(\neg A)$ in these cases as our modification.\\

{\bf Remark 1.}  Let us note the following passage in `Kaccayanagotta Sutta' (SN12.15) \cite[p. 544]{SN}:
\begin{quote}
``All exists": this is one extreme. ``All does not exists": this is the second extreme. Without veering towards either of these extremes, the Tathagata teaches the Dhamma by the middle: ``With ignorance as condition, volitional formations come to be,...
\end{quote}
 
In the Pali 
\begin{quote}
Sabbamatth\={i}'ti kho kacc\={a}na ayameko anto. Sabba\d{m} natth\={i}'ti aya\d{m} dutiyo anto. Ete te kacc\={a}na ubho ante anupagamma majjhena tath\={a}gato dhamma\d{m} deseti. 
Avijj\={a}paccay\={a} sa\d{n}kh\={a}r\={a}. \\
Sa\d{n}kh\={a}rapaccay\={a} ...
\end{quote}
Here we have the tuple of two alternatives
$$ A, \quad \neg A$$
or, in the Pali here,
$$ A, \quad \mbox{n'}A. $$

The interpretation of this tuple is clear. Let $A$ be a formula and $V$ a valuation.
The set $L$ of all formulas is divided into the subset $L_{1/2}$ of all formulas $P$ such that $V(P)=0$ and the subset $L_{2/2}$ of formulas $P$ such that $V(P)=1$. Then $A, \neg A$ are representatives of $L_{1/2}, L_{2/2}$ respectively if $V(A)=0$, while otherwise they are representatives of $L_{2/2}, L_{1/2}$ respectively. Hence our interpretation of this dilemma is: {\it If somebody denies both two alternatives, he/she intends to mean by abbreviation that any formula of any propositional calculus in which $A$ is a formula cannot describe the reality.} 

In other words, the denial of both $A$ and $\neg A$ is nothing but the denial of $A\vee(\neg A)$, which is a tautology, and it leads to the denial of all formulas in the propositional calculus considered, since, for any formula $P$, the formula $P \rightarrow A\vee(\neg A)$ is a tautology, too. Of course this argument is an intentional confusion of the object logic and the metalogic. Anyway, this dilemma may be a prototype of the catu\d{s}ko\d{t}i or tetralemma. 

Also see `Aggi-Vacchagotta Sutta' (MN72) \cite[p. 590]{MN} in the `Majjhima Nikaya', which contains both dilemmas and tetralemmas.

\section{Creator of suffering}

As the second example we take a passage from `Acela Sutta' (SN12.17) \cite[p. 546]{SN} in the `Samyutta Nikaya'. The English translation reads:

\begin{quote}
Q1: Master Gotama, is suffering created by oneself?

A1: Not so, Kassapa.

Q2: Then, Master Gotama, is suffering created by another?

A2: Not so, Kassapa.

Q3: Then, Master Gotama, is suffering created both by oneself and by another?

A3: Not so, Kassapa.

Q4: Then, Master Gotama, has suffering arisen fortuitously, being created neither by oneself nor by another?
 
A4: Not so, Kassapa.

...

Q7: Teach me about suffering, Blessed One!

A7: ``The one who acts is the one who experiences the result of the act" amounts to the eternalist statement ``suffering is created by oneself". ``The one who acts is someone other than the one who experiences the result of the act" amounts to the annihilationist statement ``suffering is created by another". Without veering towards either of these extremes, the Tathagata teaches the Dhamma by the middle: With ignorance as condition volitional formations come to be; With volitional formations...
\end{quote}

At the moment the English translation can be formulated as
\begin{equation}
A,\quad B, \quad A\wedge B, \quad \neg(A\vee B),
\end{equation}
but we modify it as 
\begin{equation}
A,  \quad B, \quad  A\vee B,\quad \neg(A\vee B)
\end{equation}
by replacing the third alternative as in \S 2. Here $A$ stands for `suffering is created by oneself " and $B$ stands for `suffering is created by another'. According to A7, we could assume that $B$ is equivalent to $\neg A$. But according to Q4, it seems that one can consider `suffering arises fortuitously (or without any cause, as a result of chance) ' even if suffering is created neither by oneself nor by another. Therefore we do not assume that $B$ is equivalent to $\neg A$. Our modification means that Q3 is replaced by\\

Q3': Then is it created either by oneself or by another? Can both be allowed?\\

This modified catu\d{s}ko\d{t}i can be characterized as follows. Let $A$ and $B$ be formulas in a propositional calculus $\Sigma$. First we assume that there are a valuation $V_0$ such that $V_0(A)=0$ and $V_0(B)=1$ and a valuation $V_1$ such that $V_1(A)=1$ and $V_1(B)=0$. If it is the case, let us say that $A$ and $B$ is {\bf separable} or {\bf independent}. Of course it is the case if $A$ and $B$ are distinct sentence letters of $\Sigma$. For any formula $P$ the pair of the truth values $(V_0(P), V_1(P))$ should be one of $(0,1), (1,0), (1,1), (0,0)$. Clearly the formulas of the tuple (7) are representatives of the four possible cases. The situation is same as in \S 2, where $B$ is $\neg A$. But we do not assume here that $A\vee B$ is identically true and $\neg(A\vee B)$ is identically false. Therefore our interpretation is that, {\it when somebody denies all the four alternatives of (7), then he/she intends to mean in abbreviation that the reality cannot be described by any formula of any propositional calculus in which $A$ and $B$ are independent formulas.}

By the way the original Pali text reads:

\begin{quote}
Q1: Kinnu kho bho gotama, saya\d{m} kata\d{m} dukkhanti?

A1: M\={a} heva\d{m} kassap\={a}.

Q2: Kimpana bho gotama, parakata\d{m} dukkhanti?

A2: M\={a} heva\d{m} kassap\={a}.

Q3: Kinnu kho bho gotama, saya\d{m} kata\~{n}ca parakata\~{n}ca dukkhanti?

A3: M\={a} heva\d{m} kassap\={a}.

Q4: Kimpana bho gotama, asaya\d{m}k\={a}ra\d{m} aparak\={a}ra\d{m} adhiccasamuppanna\d{m} dukkhanti?

A4: M\={a} heva\d{m} kassap\={a}.
\end{quote} 
Therefore the tuple of the four alternatives in Pali is
\begin{equation}
A, \quad B, \quad A\mbox{-ca}\  B\mbox{-ca}, \quad \mbox{a-}A 
\  \mbox{a-}B.
\end{equation}
So we wonder whether `$A$-ca $B$-ca' can be formalized as $ A \vee B$ or not.

The Chinese translation (`Taisho Tripitaka' 99-302) reads:
\begin{quote}
Q1: Y\'{u}nh\'{e} Q\'{u}t\'{a}n , k\v{u} z\`{i}zu\`{o} y\'{e}?

A1: K\v{u} z\`{i}zu\`{o} zhe, c\v{i} sh\`{i} w\'{u}j\`{i}.

Q2: Y\'{u}nh\'{e} Q\'{u}t\'{a}n, k\v{u} t\={a}zu\`{o} y\'{e}?

A3: K\v{u} t\={a}zu\`{o} zhe, c\v{i} y\`{i} w\'{u}j\`{i}.

Q3: K\v{u} z\`{i}t\={a}zu\`{o} y\'{e}?

A3: K\v{u} z\`{i}t\={a}zu\`{o}, c\v{i} y\`{i} w\'{u}j\`{i}.

Q4: Y\'{u}nh\'{e} Q\'{u}t\'{a}n, k\v{u} f\={e}i z\`{i} f\={e}i t\={a},
w\'{u}y\={i}n zu\`{o} y\'{e}?

A4: K\v{u} f\={e}i z\`{i} f\={e}i t\={a}, c\v{i} y\`{i} w\'{u}j\`{i}.
\end{quote}
Thus the tuple of the four alternative in Chinese is
\begin{equation}
A, \quad B, \quad A\   B,\quad \mbox{f\={e}i}\   A\  \mbox{f\={e}i}\  B.
\end{equation}
We note that no conjunctions like `and' `or' do not appear explicitly.\\

{\bf Remark 2.}  Since $A\wedge B$ entails $A\vee B$, one who denies the third alternative of our modified catu\d{s}ko\d{t}i simultaneously denies the third alternative of the usual catu\d{s}ko\d{t}i.
 (Note that $A\wedge B$ is identically false if $B=\neg A$, but otherwise it can take the truth value $1$ for a certain valuation.) According to
\cite[p. 11]{Wayman}, Tson-kha-pa's annotation to `Madhyamakak\={a}rik\={a}' I, 1 explains that $A\wedge B$ is the philosophical position of Nyaya-Vaisesika school, and $\neg(A\vee B)$ is that of Lokayata school, where $A$ stands for ``It arises from itself", $B$ stands for ``It arises from other " and ``It arises without cause (or by chance)" entails $\neg(A\vee B)$. Therefore if somebody denies all $A, B, A\vee B, \neg(A\vee B)$, he/she denies $A\wedge B$ a fortiori. 

{\bf Remark 3.} Let us note the following passage in `Samanupassana Sutta' 
(SN22.47) \cite[p. 886]{SN} in the `Samyutta Nikaya':
\begin{quote}
The thought  ``I will be percipient", ``I will be non-percipient'' and ``I will be neither percipient nor non-percipient" --these do not occur.
\end{quote}
Or in the Pali
\begin{quote}
Ayamahamasmiti'pissa na hoti, bhavissanti'pissa na hoti, na bhavissanti'pissa na hoti, sa\~{n}\~{n}\={i} bhavissanti'pissa na hoti, asa\~{n}\~{n}\={i} bhavissanti'pissa na hoti, nevasa\~{n}\~{n}\={i}n\={a}sa\~{n}\~{n}i bhavissanti'pissa na hot\={i}ti.
\end{quote} 

Here we find the trilemma
\begin{equation}
A,\quad B, \quad \neg(A\vee B),
\end{equation}
where $A$ stands for `I will be percipient' and $B$ stands for `I will be non-percipient", which formalizes 
$$A, \quad B, \quad \mbox{neither}\  A\  \mbox{nor}\  B.
$$
Or the trilemma
\begin{equation}
A,\quad \neg A, \quad \neg(A\vee(\neg A))
\end{equation}
formalizes
$$
A, \quad \mbox{a-}A,\quad \mbox{neva-}A\  \mbox{n\={a}-}A,
$$
where $A$ stands for `sa\~{n}\~{n}\={i} bhavissanti'. 

The trilemma (10) lacks the third alternative of the tetralemma (7). However, when somebody denies all the alternatives of the trilemma (10), he/she implicitly denies the third alternative $A\vee B$ of the tetralemma (7) too, since
the denial of both $A$ and $B$ entails the denial of $A\vee B$ provided that we hold the classical logic as the metalogic. Therefore the interpretation of the trilemma (10) or (11) is the same as that of the tetralemma (7) or (3). In other words this trilemma is equivalent to the tetralemma.

\section{Dual modification} 

An alternative modification of catu\d{s}ko\d{t}i could be given by replacing (7) by the tuple
\begin{equation}
A,\quad B,\quad A\wedge B,\quad \neg(A\wedge B).
\end{equation}
In this tuple (12) the third alternative coincides with the usual translation of catu\d{s}ko\d{t}i, but the fourth alternative is formalized in different way to the usual one. 

When $B=\neg A$, the tuple (12) turns out to be 
\begin{equation}
A,\quad \neg A,\quad A\wedge(\neg A),\quad \neg(A\wedge(\neg A))
\end{equation}
instead of (3). The components of the tuple (13) are representatives of $L_1, L_2, L_4, L_3$ respectively, where $L_j, j=1,2,3,4,$ are the subsets of formulas defined in \S 2. Thus the interpretation of this alternative modification of catu\d{s}ko\d{t}i is the same, that is,  {\it If 
somebody denies all these alternatives, then he/she intends to mean by abbreviation using the representatives that he/she declares that the ultimate truth or the reality cannot be described by any formula of any propositional calculus in which $A$ is a generic formula. Particularly any propositional calculus for which $A$ is a sentence letter doesn't work to describe the reality.} 

In fact the pair of the third and fourth alternatives of (13) is the mere exchange of those of (3). On the other hand, the relation between (12) and (7) can be different from a mere exchange of order if $B$ is not $\neg A$. But the tuple of the subsets of formulas considered in \S 4 represented by the components of (12) coincides with those of (7) except for the order. 

The possibility of this formalization $\neg(A\wedge(\neg A))$ of the fourth alternative of catu\d{s}ko\d{t}i was discussed by J. Westerhoff in \cite[footnote of page 375]{Westerhoff}. He discusses not on Pali texts but on Sanskrit text by N\={a}g\={a}rjuna, and his opinion seems that this interpretation is impossible. \\

Let us note that in `Sikkha Sutta' (AN 4.99) \cite[p. 479]{AN} of the `Anguttara Nikaya' the Blessed One says:
\begin{quote}
Bhikkhus, there are four kinds of persons found existing in the world. What four? One who is practicing for his own welfare but not for the welfare of others; one who is practicing for the welfare of others but not for his own welfare; one who is practicing neither for his own welfare nor for the welfare of others; and one who is practicing both for his own welfare and for the welfare of others.
\end{quote}
Here we can find the tetralemma
$$A\wedge(\neg B),\quad 
(\neg A)\wedge B, \quad 
(\neg A)\wedge(\neg B), \quad A\wedge B,$$
where $A$ stands for `he is practicing for his own welfare' and $B$ stands for `he is practicing for the welfare of others'. 
Note that in the Pali this passage reads:
\begin{quote}
Attahitt\={a}ya pa\d{t}ipanno no parahit\={a}ya; Parahit\={a}ya pa\d{t}ipanno no attahit\={a}ya; Neva attahit\={a}ya ca pa\d{t}ipanno no parahit\={a}ya; Attahit\={a}ya ca pa\d{t}ipanno parahit\={a}ya ca.
\end{quote}
That is
$$A\  \mbox{no}\  B, \quad B\  \mbox{no}\  A, \quad \mbox{neva}\  A\  \mbox{ca no}\  B,
\quad A\  \mbox{ca}\  B\  \mbox{ca}. 
$$
The Chinese translation of the similar `Valahaka Sutta' (AN 4.102) is found as Taisho 125-25.10. The tuple in the Chinese is
$$A\  \mbox{\'{e}r b\`{u}}\  B, \quad
B\  \mbox{\'{e}r b\`{u}}\  A, \quad
\mbox{y\`{i} b\`{u}}\  A\  \mbox{y\`{i} b\`{u}}\  B, \quad
\mbox{y\`{i}}\  A\  \mbox{y\`{i}}\  B.
$$

Here the Blessed One does not intend to deny the four alternatives, but intends merely to classify people, although He may agree with the opinion that one which satisfies $A\wedge B$ is the most excellent and sublime. 
It is not the case that the four alternatives are affirmed simultaneously for a single person. By the Chinese translation Taisho 125-25.10 of AN 4.102 this is explicitly expressed as 
\begin{quote}
Hu\`{o} y\v{o}u y\'{u}n l\'{e}i \'{e}r b\`{u} y\`{u}, hu\`{o} y\v{o}u y\'{u}n y\`{u} \'{e}r b\`{u} l\'{e}i, hu\`{o} ..., 
\end{quote}
where ``y\v{o}u y\'{u}n l\'{e}i \'{e}r b\`{u} y\`{u}" means `there is a cloud which thunders and does not rain' and so on, and ``hu\`{o}" means the disjunction, that is, ``Hu\`{o} ..., hu\`{o }..."
means ``Either ... or ....", or, more precisely speaking, ``On the
one hand ..., on the other hand,..." in this context. Here ``\'{e}r" = `and', ``hu\`{o}" = `or' are explicit conjunction words.

Although we cannot give a clear example in which all the alternatives are affirmed simultaneously, we would like to spend few words about affirmative catu\d{s}ko\d{t}i in which all the four alternatives are affirmed.\\  

Let the tuple (7) be called {\bf the modified catu\d{s}ko\d{t}i} generated by $A$, $B$ and the tuple (12) be called {\bf the dual modified catu\d{s}ko\d{t}i} generated by $A$, $B$. Therefore the tuple (3) is the modified catu\d{s}ko\d{t}i generated by $A, \neg A$, and (13) is the dual modified catu\d{s}ko\d{t}i generated by $A, \neg A$.

Then it is easy to see under the classical propositional calculus that the dual 
modified catu\d{s}ko\d{t}i (13) generated by $A, \neg A$ is equivalent to the modified catu\d{s}ko\d{t}i generated by $A, \neg A$ except for the exchange of the order of the alternatives, since
$$A\wedge(\neg A) \Leftrightarrow \neg(A\vee(\neg A)),\quad
\neg(A\wedge (\neg A))\Leftrightarrow A\vee(\neg A).$$

Now suppose that somebody denies all the alternatives of the dual modified catu\d{s}ko\d{t}i generated by $\neg A, \neg B$, that is,
$$ \neg A,\quad \neg B,\quad (\neg A)\wedge(\neg B), \quad \neg((\neg A)\wedge(\neg B)).$$
If we formalize this metalogical denial by the operation $\neg$ on all the alternatives in the object symbol logic, then the result is easily seen to be the modified catu\d{s}ko\d{t}i generated by $A, B$, that is, (7). In this sense the affirmation of (7) is nothing but the negation of the {\bf dual} catu\d{s}ko\d{t}i generated by $\neg A, \neg B$. Here the affirmation means affirmation of all the alternatives, and the negation means denial of all the alternatives. We note that $A, B$ are independent if and only if $\neg A, \neg B$ are independent. 

Therefore we can say that {\it the affirmation of the modified catu\d{s}ko\d{t}i (3) is nothing but the negation of (3) itself }. Of course this argument is a confusion of the object logic and the metalogic, but the conclusion, the coincidence of affirmation with negation, is a dialectical situation in a sense. So, if we want to formalize this argument as the total, we should adopt a paraconsistent logic, maybe.\\

{\bf Remark 4.} N\={a}g\={a}rjuna's M\={u}lamadhyamakak\={a}rik\={a} XVIII.8 is a problematic verse, which reads:
\begin{quote}
Everything is real and is not real, \\
Both real and not real, \\
Neither real nor not real.\\
This is Lord Buddha's teaching.
\end{quote}

According to J. L. Garfield, \cite[p. 250]{Garfield}, in contrast with \cite[p.113]{Inada},
this verse is an example of affirmative catu\d{s}ko\d{t}i without intention of denial, and N\={a}g\={a}rjuna here intends merely to mean ``Everything is {\it conventionally} real, and is {\it ultimately} unreal; Everything has both characteristics; Nothing is ultimately real". 
Therefore the opinion of Garfield may be that a reading of this verse as an affirmative catu\d{s}ko\d{t}i which is equivalent to the negative catu\d{s}ko\d{t}i dialectically as above is a nihilistic one which is very hard to sustain. See \cite[Footnote 93, p.251]{Garfield}. 

\section{Finitude and infinitude of the world}

We can find the following passage called ant\={a}natav\={a}da argument in `Brahmaj\={a}la Sutta' (DN 1) in the `Digha Nikaya' \cite[p.78]{DN}:
\begin{quote}
P0: There are some ascetics or brahmins who proclaim the finitude and infinitude of the world on four grounds. What four?

P1: A certain ascetic or brahmin thinks: ``This world is finite and bounded. ({\it  [Pali]} antav\={a} aya\d{m} loko pariva\d{t}umo; {\it [Chinese]} sh\`{i}ji\={a}n y\v{o}ubi\={a}n. )"

P2: A certain ascetic or brahmin thinks: ``This world is infinite and unbounded. ({\it [Pali]} anato aya\d{m} loko apariyanto; {\it [Chinese]} sh\`{i}ji\={a}n w\'{u}bi\={a}n.) Those who say it is finite are wrong. "

P3: A certain ascetic or brahmin, perceiving the world as {\it finite up-and-down, and infinite across} ({\it [Chinese]} sh\`{a}ngf\={a}ng y\v{o}ubi\={a}n s\`{i}f\={a}ng w\'{u}bi\={a}n ),  thinks: ``The world is finite and infinite. ({\it [Pali]} antav\={a} ca aya\d{m} loko ananto ca; {\it [Chinese]} sh\`{i}ji\={a}n y\v{o}ubi\={a}n w\'{u}bi\={a}n. ) Those who say it is finite are wrong, and those who say it is infinite are wrong."

P4: A certain ascetic or brahmin argues: ``This world is neither finite nor infinite. ({\it [Pali]} nev\={a}ya\d{m} loko antav\={a} na pan\={a}nanto; {\it [Chinese]} sh\`{i}ji\={a}n f\={e}i y\v{o}ubi\={a}n f\={e}i w\'{u}bi\={a}n.) Those who say it is finite are wrong, and so those who say it is infinite, and those who say it is finite and infinite."

P5: These are the four ways. There is no other way.
\end{quote} 

Let us try to formalize this argument.\\

Consider the tuple
\begin{equation}
A\wedge(\neg B),\quad
(\neg A)\wedge B, \quad
A\wedge B, \quad
(\neg A)\wedge(\neg B),
\end{equation}
where $A$ and $B$ are formulas. 
(This tuple is that of AN 4.99 mentioned in \S 5 except for the exchange of the order of the alternatives.) We can consider that this is {\bf the proper (unmodified) catu\d{s}ko\d{t}i}. Hereafter we denote by $C_1, C_2, C_3, C_4$ the alternatives of the tuple (14). It is easy to verify the properties
\begin{equation}
C_i\wedge C_j \Leftrightarrow \curlywedge \quad \mbox{if}\  i\not= j,
\end{equation}
and
\begin{equation}
C_1\vee C_2\vee C_3\vee C_4 \Leftrightarrow \curlyvee,
\end{equation}
by using the auxiliary truth-value table:
$$
\begin{array}{ll|llll}
A & B & C_1 & C_2 & C_3 & C_4 \\ \hline
0 & 0 & 0 & 0 & 0 & 1 \\
0 & 1 & 0 & 1 & 0 & 0 \\
1 & 0 & 1 & 0 & 0 & 0 \\
1 & 1 & 0 & 0 & 1 & 0
\end{array}
$$

Of course, when $B=\neg A$, the tuple (14) is equivalent to (2), 
say, $A, \neg A, \curlywedge, \curlywedge$, and this is not interesting as pointed in \S 1. However we can consider the case in which $A$ is $\exists xFx$ and $B$ is $\exists x\neg Fx$, where $F$ is a one-place predicate and $x$ is a variable. Here we adopt the classical predicate calculus developed in \cite[Chapter 4]{Sider}. In this case it is easy to verify that the tuple (14) is equivalent to
\begin{equation}
\forall xFx, \quad
\forall x\neg Fx, \quad
(\exists xFx)\wedge (\exists x\neg Fx), \quad
\forall x(Fx\wedge(\neg Fx)),
\end{equation}
since 
$$\neg\exists x\neg Fx \Leftrightarrow
\forall xFx,\quad
\neg\exists xFx\Leftrightarrow \forall x\neg Fx.
$$

Here actually we have $C_4\Leftrightarrow \curlywedge$ but $C_3$ can be nonequivalent to $\curlywedge$ when the domain of the variable  $x$ contains distinct elements. 

So, this catu\d{s}ko\d{t}i may formalize the ant\={a}natav\={a}da argument very well, if we consider that $Fx$ stands for `the world is finite and bounded with respect to the direction $x$'. In fact, if $a$ stands for `up-and-down' and $b$ stands for `east-west-south-and-north', the third ascetic or brahmin believes that both $Fa$ and $\neg Fb$ are true, therefore, $C_3$ is true. Moreover we note that the tuple $$C_1,\quad C_2,\quad C_3,\quad C_4
$$ is clearly
equivalent to the tuple
$$C_1, \quad
C_2\wedge(\neg C_1), \quad
C_3\wedge(\neg C_1)\wedge(\neg C_2), \quad
C_4\wedge(\neg C_1)\wedge(\neg C_2)\wedge(\neg C_3)
$$
as described in the text of the Sutta little bit redundantly.

In view of (15)(16) the saying P5 of the Blessed One is exact. If somebody denies all the alternatives $C_i, i=1,2,3,4,$ as not to be attached, then the result is the absolute empty, or `{\it nibbuti}' (={\it nibb\={a}na}, perfect peace beyond reasoning) and  `{\it anup\={a}d\={a}-vimutta}' (emancipation without clinging). \\ \medskip

{\bf\large Acknowledgment }\  
The author would like to express his sincere thanks to Professor Yasuo Deguchi (Kyoto University), who introduced the study subject `catu\d{s}ko\d{t}i' to the author
and encouraged him during the preparation of this note.
But it is regrettable that the author has not been able to 
make the most of his comments, which are philosophically profound.

The author would like to express his sincere thanks to Doctor
Takuro Onishi for giving helpful comments to
ameliorate the manuscript through the occasion
of a discussion meeting, which was financially supported by the Department of Philosophy, Kyoto University. \\ \medskip

{\bf\Large Appendix}\\

Let us consider a propositional calculus $\Sigma$ and an arbitrary 
pair of valuations 
$v_1$ and $v_2$ in $\Sigma$. 
Let us denote $\mathfrak{v}(P)=(v_1(P), v_2(P))$ for any formula $P$ of
$\Sigma$.
 Then $\mathfrak{v}(P)$ can take one of the four vector values $(1,0), (0,1), (1,1), (0,0)$. Let us denote
$$\mathfrak{b}=(1,0),\quad \mathfrak{n}=(0,1),
\quad \mathfrak{t}=(1,1), \quad
\mathfrak{f}=(0,0).$$

Now, by tedious calculations,  it can be verified that the performance of
the four truth values $\mathfrak{v}$ obeys the following tables:

$$
\begin{array}{r|r}
\neg & \\ \hline
\mathfrak{t} & \mathfrak{f} \\
\mathfrak{b} & \mathfrak{n} \\
\mathfrak{n} & \mathfrak{b} \\
\mathfrak{f} & \mathfrak{t}
\end{array}
\qquad
\begin{array}{r|rrrr}
\vee & \mathfrak{t}&\mathfrak{b}&\mathfrak{n}&\mathfrak{f}\\ \hline
\mathfrak{t}&\mathfrak{t}&\mathfrak{t}
&\mathfrak{t}&\mathfrak{t} \\
\mathfrak{b}&\mathfrak{t}&\mathfrak{b}
&\mathfrak{t}&\mathfrak{b} \\
\mathfrak{n}&\mathfrak{t}&\mathfrak{t}
&\mathfrak{n}&\mathfrak{n}\\
\mathfrak{f}&\mathfrak{t}&\mathfrak{b}
&\mathfrak{n}&\mathfrak{f}
\end{array} 
$$
$$
\begin{array}{r|rrrr}
\wedge & \mathfrak{t}&\mathfrak{b}&\mathfrak{n}&\mathfrak{f}\\ \hline
\mathfrak{t}&\mathfrak{t}&\mathfrak{b}
&\mathfrak{n}&\mathfrak{f} \\
\mathfrak{b}&\mathfrak{b}&\mathfrak{b}
&\mathfrak{f}&\mathfrak{f} \\
\mathfrak{n}&\mathfrak{n}&\mathfrak{f}
&\mathfrak{n}&\mathfrak{f}\\
\mathfrak{f}&\mathfrak{f}&\mathfrak{f}
&\mathfrak{f}&\mathfrak{f}
\end{array} 
\qquad
\begin{array}{r|rrrr}
\rightarrow & \mathfrak{t}&\mathfrak{b}&\mathfrak{n}&\mathfrak{f}\\ \hline
\mathfrak{t}&\mathfrak{t}&\mathfrak{b}
&\mathfrak{n}&\mathfrak{f} \\
\mathfrak{b}&\mathfrak{t}&\mathfrak{t}
&\mathfrak{f}&\mathfrak{f} \\
\mathfrak{n}&\mathfrak{t}&\mathfrak{b}
&\mathfrak{t}&\mathfrak{b}\\
\mathfrak{f}&\mathfrak{t}&\mathfrak{t}
&\mathfrak{t}&\mathfrak{t}
\end{array} 
$$

Therefore we have a semantics for a four-valued logic, which is 
similar to that of Dunn for \textbf{FDE} (Finite Degree Entailment) adopted by G. Priest \cite[\S 3.2]{Priest} with the Hasse diagram in \cite[p.33]{Priest}.
The different point is that the value of a negation $\neg$ is fixed for $\mathfrak{b}$ and
$\mathfrak{n}$ while it toggles $\mathfrak{t}$ and 
$\mathfrak{f}$ in the \textbf{FDE} semantics but it toggles $\mathfrak{b}$ and $\mathfrak{n}$,
too, for our case. On the other hand, suggested by Y. Deguchi, J. Garfield and G. Priest, 2008, \cite[pp. 398-399]{DeguchiGP}, A. J. Cotnoir \cite{Cotnoir}
introduces the semantics \textbf{B4},
in which the value of a negation $\neg$ toggles $\mathfrak{b}$ and
$\mathfrak{n}$, too. In other words,
using the symbol and the interpretations of \cite{Cotnoir}, 
we can put
\begin{align*}
&\langle 1,1\rangle (=\mbox{both CT and UF}):=\mathfrak{b}=(1,0) \\
&\langle1,0\rangle(=\mbox{CT but not UF}):=\mathfrak{t}=(1,1) \\
&\langle0,1\rangle (=\mbox{not CT but UF}):=\mathfrak{f}=(0,0) \\
&\langle0,0\rangle (=\mbox{neither CT nor UF}):=\mathfrak{n}=(0,1),
\end{align*}
where `CT' stands for `conventionally true' and `UF' stands for
`ultimately false'. Then the semantics for our pairing of valuations $\mathfrak{v}=(v_1,v_2)$ coincides with
that of \textbf{B4}. 
(Note that if, we are not sure but, `not UF' is equivalent to `ultimately true',
then 
$v_1(A)=1\  [\!( 0 )\!]$ iff the proposition $A$ is {\bf conventionally true} [\!( false )\!] and $v_2(A)=1\  [\!( 0 )\!] $ iff 
the proposition $A$ is {\bf ultimately true} [\!( false )\!].) In this sense
our interpretation of the modified catu\d{s}ko\d{t}i,
if we take $v_1=V_0, v_2=V_1$, 
is compatible with the paraconsistent point of view of 
Y. Deguchi, J. Garfield and G. Priest formalized by A. J. Cotnoir.
 
But a more appropriate semantical formulation of catu\d{s}ko\d{t}i
has been presented by T. Onishi,
who appeals to the concept `bilattice'. The details will be given in \cite{Onishi}.

\end{document}